\documentclass{amsart}
\usepackage{amssymb,amsfonts,amsmath,amsthm}
\usepackage[all]{xy}
\usepackage[shortlabels]{enumitem}

\newtheorem{theorem}{Theorem}

\newtheorem{corollary}[theorem]{Corollary}

\newtheorem*{theoremA}{Theorem A}
\newtheorem*{corollaryA}{Corollary A}
\newtheorem*{theoremB}{Theorem B}
\newtheorem*{theoremC}{Theorem C}

\newcommand{\real}{\mathbb{R}}

\newcommand{\si}{\sigma}

\newcommand{\la}{\lambda}

\newcommand{\vp}{\varphi}
\newcommand{\om}{\omega}
\newcommand{\Om}{\Omega}
\newcommand{\na}{\nabla}

\newcommand{\de}{\delta}
\newcommand{\De}{\Delta}
\newcommand{\ka}{\kappa}
\newcommand{\ti}{\tilde}
\newcommand{\vol}{\text{\rm vol}}

\newcommand{\p}{\partial}

\newcommand{\m}{\mathcal}

\newcommand{\ric}{\mathrm{Ric}}
\newcommand{\dv}{\mathrm{div}}

\title[]{Rigidity Theorems of conformal class on compact manifolds with boundary} 
\author{Ezequiel Barbosa} \address{Instituto de Ci\^encias Exatas, Universidade Federal de Minas Gerais, CEP 31270-901, Belo Horizonte, MG, Brasil}\curraddr{Imperial College London, Huxley Building, 180 QueenÕs Gate, London SW7 2RH, United Kingdom} \email{ezequiel@mat.ufmg.br} \author{Heudson Mirandola} \address{Instituto de Matematica, Universidade Federal do Rio de Janeiro, CEP 21945-970, Rio de Janeiro, RJ, Brazil}\curraddr{Imperial College London, Huxley Building, 180 QueenÕs Gate, London SW7 2RH, United Kingdom}
\email{mirandola@im.ufrj.br}
\author{Feliciano Vitorio}\address{Instituto de Matem\'atica, Universidade Federal de Alagoas, CEP 57072-900, Macei\'o, Al, Brazil}\email{feliciano@pos.mat.ufal.br}

\thanks{This work is partially supported by CNPq.}
\subjclass[2000]{Primary ; Secondary }

\begin{document}
\maketitle
\begin{abstract} Let $M$ be a compact manifold with boundary. In this paper, we discuss some rigidity theorems of metrics in a same conformal class that fixes the boundary and satisfy certain integral conditions on the the scalar curvatures and the mean curvatures on the boundary. No condition on the first eigenvalues of operators is need. 
\end{abstract}
\section{Introduction}

Let $(M, g_0)$ be a compact $n$-dimensional Riemannian smooth manifold with $n\ge 2$ and  nonempty smooth boundary $\p M$ (possibly non connected). Let $R_{g_0}$ denote the scalar curvature of $(M,g_0)$ and let $H_{g_0}=\dv_{g_0} \eta_{g_0}$ denote the mean curvature of $\p M$ in $(M,g_0)$, in the direction of the exterior conormal $\eta=\eta_{g_0}$. If  $n=2$ then $K_{g_0} = R_{g_0}/2$ denotes the Gaussian curvature and $H_{g_0}=\ka_{g_0}$ denotes the geodesic curvature of $\p M$ with respect to $g_0$. 

Recall the conformal class of a metric $g$ on $M$,  say $[g]$, is the set of metrics of the form $\ti g=\mu g$, where $\mu$ is a positive smooth function defined on $M$. Escobar \cite{E1} dealt with the following question:  
\begin{quotation} {\it Given a metric $g\in [g_0]$ with $R_{g}=R_{g_0}$ in $M$, and $H_{g}=H_{g_0}$ on $\p M$, when is $g=g_0$?}
\end{quotation} 
As Escobar observed in \cite{E2}, this question does not have a positive answer in general. Indeed, he gave the description of the conformally flat metrics $g\in [\de_{ij}]$ on the ball $B=\{x\in \real^n \mid |x|\le 1\}$, with $n\ge 3$, having constant scalar curvature and constant mean curvature on $\p B$. By this classification theorem, there is a non compact set of metrics $g\in [\de_{ij}]$ with $R_g=0$ and $H_g=1$. Also in \cite{E2}, Escobar showed that on the annulus $A_{a,b}=\{x\in \real^n \mid a<|x|<b\}$ there are several metrics $g\in [\de_{ij}]$ with the same constant scalar curvature and the same constant mean curvature on the boundary, provided $a/b$ is big enough.

In terms of PDE's,  the scalar curvatures and mean curvatures of the metrics $g_0$ and $g\in [g_0]$ satisfy the following differential equations. For $n\ge 3$, we write $g=u^{\frac{4}{n-2}}g_0$. It holds: 
\begin{equation}\label{equations>=3}
\left\{
\begin{array}{l}
L_{g_0}u + c(n) R_g u^{\frac{n+2}{n-2}}=0, \mbox{ in } M,\\
B_{g_0}u - 2c(n)u^\frac{n}{n-2}=0, \mbox{ on } \p M,
\end{array}\right.
\end{equation} 
where $c(n)=\frac{n-2}{4(n-1)}$, $L_{g_0}(u)=\De_{g_0}u -c(n) R_{g_0}u$ and $B_{g_0}(u)=\frac{\p u}{\p \eta} + 2c(n) H_{g_0} u
$.  For $n=2$, we write $g=e^{2u}g_0$. It holds: 
\begin{equation}\label{equations=2}
\left\{
\begin{array}{l}
L_{g_0}u + K_g e^{2u}=0, \mbox{ in } M,\\
B_{g_0}u-k_g e^u=0, \mbox{ on } \p M,
\end{array}\right.
\end{equation} 
where $L_{g_0}(u)=\De_{g_0}u -K_{g_0}$ and $B_{g_0}(u)=\frac{\p u}{\p \eta} + \ka_{g_0}$. 
The second pair of operators Escobar talked about are defined by
\begin{equation}\label{equation-linear}
\left\{
\begin{array}{l}
\m L_{g_0} = \De_{g_0} + \frac{1}{n-1}R_{g_0}, \mbox{ in } M,\\
\m B_{g_0} = \frac{\p}{\p\eta_{g_0}} - \frac{1}{n-1}H_{g_0}, \mbox{ on } \p M.
\end{array}\right.
\end{equation} 
The operators $(\m L_{g_0},\m B_{g_0})$ are the linearizations of (\ref{equations>=3}), when $n\ge 3$, at $u=1$, and (\ref{equations=2}), when $n=2$, at $u=0$, for the cases $R_g=R_{g_0}$ and $H_g=H_{g_0}$. Escobar \cite{E1} proved the following
\begin{theoremA}[Theorem 1 of \cite{E1}] Let $g\in [g_0]$ with $R_g=R_{g_0}$ and $H_g=H_{g_0}\le 0$. If both first eigenvalues $\la_1(\m L_{g},\m B_{g})$ and $\la_1(\m L_{g_0},\m B_{g_0})$ are positive or one of them is equal to zero then $g=g_0$.
\end{theoremA}

As a consequence, he obtained the following
\begin{corollaryA}[Corollary 2 of \cite{E1}] Let $g\in [g_0]$ satisfying $R_g=R_{g_0}\le 0$ and $H_g=H_{g_0}\le 0$. Then $g=g_0$.
\end{corollaryA}

Min-Oo \cite{Min}  conjectured the following: {\it Let $(M,g)$ be an $n$-dimensional compact Riemannian manifold with boundary and  scalar curvature $R_g\ge n(n-1)$. Assume the boundary is  isometric to the standard sphere $S^{n-1}$ and is totally geodesic in $M$. Then, $(M,g)$ is isometric to the upper hemisphere $S^n_+$.} Min-Oo's conjecture fell in 2011, when Brendle, Marques and Neves exhibited a beautiful  counterexample. On the other hand, Hang and Wang \cite{HW1} proved Min-Oo's conjecture is true for the case the metric is conformal to the metric of $S_+^n$. They proved the following 

\begin{theoremB}[Theorem 3.4 of \cite{HW1}] Let $g\in [g_{S^n}]$ on $S^n_+$.  Assume that the scalar curvature $R_g\ge R_{g_{S^n}}=n(n-1)$ and $g=g_{S^n}$ on the boundary $\p S^n_+$. Then $g=g_{S^n}$.
\end{theoremB}

Based on Theorem B and Min-Oo's Problem, we are interested into the following question: \begin{quotation} {\it Given a metric $g\in [g_0]$ with $g=g_0$ on $\p M$ and $R_g\ge R_{g_0}$ in $M$, when is $g=g_0$?}\end{quotation}

The upper hemisphere $S^n_+$ is a static manifold, that means there is a smooth function $f$ satisfying the equation
\begin{equation} \label{eq-static}
\left\{
\begin{array}{l}
f\ric - \na^2 f + (\De f) g=0, \mbox{ in } M\setminus \p M\\
f>0 \mbox{ in } M\setminus \p M, \mbox{ and } f=0, \mbox{ on } \p M
\end{array} \right.
\end{equation}
As a solution of (\ref{eq-static}) for $S^n_+$, we take, for instance, the height function $f(x)=x_{n+1}$, for all $x=(x_1,\ldots,x_{n+1})\in S^n_+$. By taking the trace in (\ref{eq-static}), we see  static manifolds are solutions of $\m L_g f=0$, for some $f\in C^2(M)$ that is positive in $M$ and vanishes on $\p M$.

Our first theorem is
\begin{theorem}\label{min-oo conformal} Let  $g=\mu^2 g_0$ a metric in the class $[g_0]$ such that $g=g_0$ on $\p M$. Let $f\in C^1(M)\cap C^2(M\setminus \p M)$, positive a.e. such that
\begin{equation}\label{hyp-int}
\int_M f(R_g-R_{g_0}) d\vol_{g_0} + 2\int_{\p M} f(H_g - H_{g_0}) d\m H^{n-1}_{g_0}\ge 0.
\end{equation}
If  $\int_M \m L_{g_0}f (1-\mu^{-2})d\vol_{g_0}\ge 0$ then  $g=g_0$.
\end{theorem}

A direct application of Theorem \ref{min-oo conformal}, by using $\m L_{g_0}(1)=\frac{1}{n-1}R_{g_0}$, is 

\begin{corollary}\label{cor f=1} Let $g=\mu^2g_0$ be a metric in the conformal class $[g_0]$ such that $g=g_0$ on $\p M$. Assume that 
\begin{equation}\label{eq cor 1}
\int_M (R_g-R_{g_0}) d\vol_{g_0} + 2\int_{\p M} (H_g - H_{g_0}) d\m H^{n-1}_{g_0}\ge 0.
\end{equation}
If $\int_M R_{g_0}(1-\mu^{-2})d\vol_{g_0}\ge 0$ then $g=g_0$.
\end{corollary}

For static metrics, it holds $\m L_{g_0}f=0$, for some $f\in C^2(M)$ that is positive in $M\setminus \p M$ and  vanishes along $\p M$. Thus, we have 

\begin{corollary}\label{corol-static} Let $g_0$ be a static metric on $M$ and $g\in [g_0]$ such that $g=g_0$ on $\p M$. If $R\ge R_0$ then $g=g_0$. 
\end{corollary}

Araujo \cite{A} studied the functional 
\begin{equation}\label{araujo}
F(g) = \int_M R_g \,d\vol_g + 2\int_{\p M} H_g \,d\si_g.
\end{equation}
restricted to the subset of metrics $\m M_{ab}=\{g\mid a\,\vol_g(M)+b\m A_g(\p M)=1\}$. He proved the critical points of $F$ are the Einstein metrics with umbilical boundary that satisfy $b(n-1)R_g=2n a H_g$.
It is worthwhile to point out assumption (\ref{eq cor 1}) of Corollary \ref{cor f=1} does not imply $F(g)\ge F(g_0)$, since volume and area elements in (\ref{eq cor 1}) does not vary with the metric. Using Gauss-Bonnet Theorem, Corollary \ref{cor f=1} in dimension 2 can be written as
\begin{corollary}\label{for f=1 n=2} Let $(M,g_0)$ be a Riemannian surface with smooth boundary. Let $u\in C^2(M)$ with $u=0$ on $\p M$ and consider the metric $g=e^{2u}g_0$. Assume 
\begin{equation*}
\int_M K_g \,d\vol_{g_0} + \int_{\p M} \ka_{g} d \m H^1_{g_0}\ge 2\pi \chi(M),
\end{equation*}
If $\int_M K_{g_0}(1-e^{-2u})d\vol_{g_0}\ge 0$ then $u=0$ in $M$.
\end{corollary}

Theorem \ref{min-oo conformal} requires no condition on the first eigenvalue of $\m L_{g_0}$ (and even $\m B_{g_0}$). However, the first eigenvalue $$\la_1(\m L_{g_0})=\inf\{\int_M(|\na\vp|^2-\frac{R_{g_0}}{n-1}\vp^2)d\vol_{g_0}\mid \vp\in C^\infty_0(M), \int_M \vp^2 d\vol_{g_0}=1\}$$  satisfies $\m L_{g_0}f+\la_1(\m L_{g_0}) f=0$, for some $C^2$ eigenfunction $f$ that is positive in $M$ and vanishes on $\p M$. Thus, Corollary \ref{cor-min-oo conformal} below is a generalization of Corollary \ref{corol-static}.  

\begin{corollary}\label{cor-min-oo conformal} Let $g\in [g_0]$  satisfying $g=g_0$, on $\p M$, and $R_g\ge R_{g_0}$. Assume   
$\la_1(\m L_{g_0})(g-g_0)\le 0$. Then $g=g_0$. \end{corollary}

Llarul \cite{Ll}, confirming a conjecture due Gromov, proved the following result:  {\it If $g$ is any metric on the sphere $S^n$  satisfying $g\ge g_0$ and $R_g\ge R_{g_{S^n}}=n(n-1)$ then  $g= g_{S^n}$}. Furthermore, for domains $\Om\subset S^n$, Hang and Wang \cite{HW2} proved the following

\begin{theoremC}[Proposition 1 of \cite{HW2}] Let $\Om$ be a smooth domain in $S^n_+$ and let $g\in [g_{S^n}]$ in $\bar \Om$, satisfying $R_g\ge n(n-1)$ and $g=g_{S^n}$ on $\p \Om$. Then either $g=g_{S^n_+}$, in $\Om$, or $g>g_{S^n_+}$ and $H<H_{g_{S^n_+}}$.   
\end{theoremC}

In this paper we prove 
\begin{theorem}\label{thm-hang-wang} Assume $R_{g_0}\ge 0$ and $\m L_{g_0}f\le 0$, for some $f\in C^2(M\setminus \p M)\cap C^1(M)$ positive a.e.. Let $\Om$ be a smooth domain in $M$ and let $g=\mu^2 g_0$, where $\mu\in C^2(\Om)\cap C^0(\bar\Om)$ is positive with $\mu|_{\p\Om}=1$. Assume that $\chi_{\{\mu<1\}}R_{g}\ge \chi_{\{\mu<1\}}R_{g_0}$. Then, it holds 
\begin{equation}\label{inequalies}
g\ge g_0 \mbox{ in } \Om, \mbox{ and } H_g\le H_{g_0} \mbox{ in } \p\Om.
\end{equation}
In addition, if $R_g\ge R_{g_0}$ then the inequalities in (\ref{inequalies}) are strict, unless $g=g_0$.
\end{theorem}

If $g_0$ is a static metric on a manifold with boundary then it holds 

\begin{corollary}\label{cor-thm-han-wan}  Let $g_0$ be a static  metric on $M$  with $R_{g_0}\ge 0$. Let $\Om\subset M$ be a smooth domain and let $g=\mu^2 g_0$, where $\mu\in C^2(\Om)\cap C^0(\bar\Om)$ is positive with $\mu|_{\p\Om}=1$. Assume that $\chi_{\{\mu<1\}}R_{g}\ge \chi_{\{\mu<1\}}R_{g_0}$. Then, it holds 
\begin{equation*}
g\ge g_0 \mbox{ in } \Om, \mbox{ and } H_g\le H_{g_0} \mbox{ in } \p\Om.
\end{equation*}
In addition, if $R_g\ge R_{g_0}$ then the inequalities above are strict, unless $g=g_0$.
\end{corollary}

As another application, by Theorem \ref{thm-hang-wang} with $f=1$, we obtain
\begin{corollary} Let $g_0$ be a metric on $M$ with $R_{g_0}=0$. Let $\Om\subset M$ be a smooth domain. Consider the metric $g=\mu^2 g_0$, where $\mu\in C^2(\Om)\cap C^0(\bar\Om)$ is positive with $\mu|_{\p\Om}=1$. Assume $\chi_{\{\mu<1\}}R_g\ge 0$. Then, it holds
\begin{equation*}
 g\ge g_0 \mbox{ in } \Om, \mbox{ and } H_g\le H_{g_0} \mbox{ in } \p\Om.
\end{equation*}
 In addition, if $R_g\ge 0$ in $\Om$, then the inequalities above are strict, unless $g=g_0$.
\end{corollary}

\section*{Aknowledgement}
The authors are very grateful to Professor Andre Neves for your suggestions and comments, and Imperial College London for the kind hospitality. During this work, E. Barbosa and H. Mirandola were supported by a CNPq Postdoctoral Fellowship, and F. Vitorio was supported by a CNPq Universal Grant.

\section{Proof of Theorem \ref{min-oo conformal}}

First, consider the case $n=2$ and write $g=e^{2u}g_0$, with  $u\in C^2(M\setminus \p M)\cap C^1(M)$. Since $g=g_0$ on $\p M$, one has  $u|_{\p M}= 0$. The geodesic curvatures $\ka_g$, $\ka_{g_0}$ satisfy
\begin{equation}\label{meancurvaturen=2}
\frac{\p u}{\p \eta} = \ka_g e^u - \ka_{g_0} = \ka_g - \ka_{g_0}, \ \mbox{ on } \p M,
\end{equation}
where $\eta = \eta_{g_0}$ is the outward unit normal vector of $(\p M,g_0)$. 
Furthermore, the Gaussian curvatures $K_g$, $K_{g_0}$ of $(M,g_0)$ and $(M,g)$, respectively, satisfy
\begin{equation}\label{conformaln=2}
\De_{g_0} u - K_{g_0} + K_{g} e^{2u}=0, \ \mbox{ in } M.
\end{equation}
Using that $u|_{\p M}=0$, by (\ref{meancurvaturen=2}) and integration by parts, we obtain 
\begin{eqnarray*}
\int_M e^{-2u} \De_{g_0} f &=& \int_M f \De_{g_0}(e^{-2u}) + \int_{\p M} (e^{-2u} \frac{\p f}{\p \eta} - f \frac{\p (e^{-2u})}{\p \eta})\\
&=& \int_M f \De_{g_0}(e^{-2u}) + \int_{\p M}\frac{\p f}{\p \eta} + 2 \int _{\p M} f \frac{\p u}{\p \eta}\\
&=& \int_M [-2f e^{-2u} (\De_{g_0}u - 2|Du|_{g_0}^2) + \De_{g_0} f] + 2 \int _{\p M} f (\ka_g-\ka_{g_0})\\
&=& \int_M [-2f e^{-2u} (K_{g_0}-K_g e^{2u} - 2|Du|_{g_0}^2) + \De_{g_0} f] 
\\&&+ \int _{\p M} 2f (\ka_g-\ka_{g_0}).
\end{eqnarray*}
Thus, since $\De_{g_0}f = \m L_{g_0}f - 2 K_{g_0}f$, we obtain
\begin{eqnarray*}
\int_M e^{-2u}(\De_{g_0}f + 2K_{g_0}f) &=& \int_M (\m L_{g_0}f + 4f e^{-2u}|Du|_{g_0}^2) 
\\&&+ \int_M 2f (K_g - K_{g_0}) +  \int _{\p M} 2f (\ka_g-\ka_{g_0}).
\end{eqnarray*}
Hence, 
\begin{eqnarray*}
\int_M \m L_{g_0}f (e^{-2u}-1) =  \int_M 4f e^{-2u}|Du|_{g_0}^2+ \int_M 2f (K_g - K_{g_0}) +  \int _{\p M} 2f (\ka_g-\ka_{g_0}).
\end{eqnarray*}

By hypothesis, $\int_M \m L_{g_0}f (1-e^{-2u})\ge 0$ and $\int_M 2f (K_g - K_{g_0}) +  \int _{\p M} 2f (\ka_g-\ka_{g_0})$. Hence, $Du=0$, which together the fact that $u|_{\p M}=0$, imply that $g=g_0$.  

Now, we assume  $n\ge 3$ and write $g=u^{\frac{4}{n-2}}g_0$, for some $u\in C^2(M\setminus \p M)\cap C^1(M)$, positive in $M$ and with $u=1$ on $\p M$. The mean curvatures $H_{g_0}=\dv_{g_0}\eta_{g_0}$ and $H_{g}=\dv_g\eta_g$ satisfy 
\begin{equation}\label{meancurvature n>2}
\frac{\p u}{\p \eta}= \frac{n-2}{2(n-1)}(H_g u^{\frac{n}{n-2}} - H_{g_0}) = \frac{n-2}{2(n-1)}(H_g - H_{g_0}), \mbox{ on } \p M,
\end{equation}
where $\eta=\eta_{g_0}$. Furthermore, the scalar curvatures $R_g$ and $R_{g_0}$ satisfy
\begin{equation}\label{conformal n>2}
\De_{g_0} u - \frac{n-2}{4(n-1)}R_{g_0}u + \frac{n-2}{4(n-1)}R_{g} u^{\frac{n+2}{n-2}}=0, \mbox{ in } M.
\end{equation}
\\
Let $\la$ be a constant to be chosen later. Using that $u|_{\p M}=1$, integrating by parts we obtain 
\begin{eqnarray*}
\int_M u^\la \De_{g_0} f &=& \int_M  f \De_{g_0}u^\la + \int_{\p M} (u^\la \frac{\p f}{\p \eta} - f \frac{\p (u^\la)}{\p \eta})\\
&=& \int_M (f \la u^{\la -1} \De_{g_0}u + f \la (\la-1) u^{\la-2} |Du|_{g_0}^2) + \int_{\p M} \frac{\p f}{\p \eta} - \la f  \frac{\p u}{\p \eta}\\
&=& \int_M [(f \la u^{\la -1}\De_{g_0}u + f\la(\la-1) u^{\la-2} |Du|_{g_0}^2) + \De_{g_0} f]
\\&& - \frac{(n-2)}{2(n-1)}\la\int_{\p M} f(H_g - H_{g_0})
\\&=& \int_M [f \la u^{\la -1} (\frac{n-2}{4(n-1)}(R_{g_0}u - R_g u^\frac{n+2}{n-2}) + f\la(\la-1) u^{\la-2} |Du|_{g_0}^2)] \\&& + \int_M \De_{g_0} f - \frac{(n-2)}{2(n-1)}\la\int_{\p M} f(H_g - H_{g_0}).
\end{eqnarray*}
Hence,
\begin{eqnarray*}
\int_M u^\la (\De_{g_0}f - \la f \frac{n-2}{4(n-1)}R_{g_0}) &=& -\la\frac{n-2}{4(n-1)} \int_M f R_g u^{\la-1+\frac{n+2}{n-2}}\\&& + \int_M \De_{g_0}f + \la(\la-1)\int_M f u^{\la-2}|Du|_{g_0}^2 \\&& - \frac{(n-2)}{2(n-1)}\la\int_{\p M} f(H_g - H_{g_0})  
\end{eqnarray*}
We choose $\la=1-\frac{n+2}{n-2}=\frac{-4}{n-2}$. We obtain
\begin{eqnarray*}
\int_M u^{\frac{-4}{n-2}} \m L_{g_0}f &=& \int_M f \frac{R_g}{n-1} + \int_M (\m L_{g_0}f - \frac{R_{g_0}}{n-1}f)
\\&&+ \frac{4(n+2)}{(n-2)^2}\int_M fu^{\frac{-2n}{n-2}}|Du|_{g_0}^2 + \frac{2}{n-1}\int_{\p M} f(H_g-H_{g_0}).
\end{eqnarray*}
We obtain
\begin{eqnarray*}
\int_M \m L_{g_0}f (u^{\frac{-4}{n-2}}-1) &=& \frac{4(n+2)}{(n-2)^2}\int_M fu^{\frac{-2n}{n-2}}|Du|_{g_0}^2 
\\&& + \frac{1}{n-1} \int_M f (R_g-R_{g_0}) + \frac{2}{n-1}\int_{\p M} f(H_g-H_{g_0}).
\end{eqnarray*}
By hypothesis, $g=u^{\frac{4}{n-2}}g_0$ satisfies $\int_M \m L_{g_0}f (1-u^{\frac{-4}{n-2}})\ge 0$. Using (\ref{hyp-int}), we obtain that $Du=0$. Since $u|_{\p M}=1$ one has $u=1$ in $M$; hence $g=g_0$. Theorem \ref{min-oo conformal} is proved.

\section{proof of Theorem \ref{thm-hang-wang}}

First, consider the case $n=2$ and write $g=e^{2u}g_0$ with $u\in C^{2}(\Om)\cap C^0(\bar \Om)$. Since $u|_{\p \Om}= 0$ the geodesic curvatures $\ka_g$, $\ka_{g_0}$ satisfy
\begin{equation}\label{meancurvature n=2}
\frac{\p u}{\p \eta} = \ka_g e^u - \ka_{g_0} = \ka_g - \ka_{g_0}, \mbox{ in } \p\Om,
\end{equation}
where $\eta = \eta_{g_0}$ is the outward unit normal vector of $(\p \Om,g_0)$. 
Furthermore, the Gaussian curvatures $K_g$, $K_{g_0}$ of $(\Om,g_0)$ and $(\Om,g)$, respectively, satisfy
\begin{equation}\label{supharmonic}
\De_{g_0} u =K_{g_0} - K_{g} e^{2u}.
\end{equation}
Let $\bar u=\min\{u,0\}$. It turns that  $\bar u$ is continuous and, in the sense of distributions, it holds $\De_{g_0}\bar u \le \chi_{\{u<0\}}\De_{g_0}u = \chi_{\{u<0\}}(K_{g_0}-K_g e^{2u})\le  \chi_{\{u<0\}}K_{g_0}(1- e^{2u})$, since $ \chi_{\{u<0\}} K_{g}\ge  \chi_{\{u<0\}} K_{g_0}$. Let $A_u=\chi_{\{u<0\}} K_{g_0}(1-e^{2u})$. We have  that $A_u=A_{\bar u}$ is a nonnegative continuous function, and
$$\De_{g_0}\bar u \le A_{\bar u}, \mbox{ in } M,$$
in the sense of distributions.  The function $A_{\bar u}$ is Lipschitz in $\bar \Om$. In fact, given $x,x\in \bar\Om$, if either $x,x_0\in \{u<1\}$, or $x,x_0\in \{u\ge 1\}$, we have $|A_{\bar u}(x)-A_{\bar u}(x_0)|=\chi_{\{u<0\}}|K_0(x)e^{2u(x)}- K_0(x_0) e^{2u(x_0)}|\le M d_{g_0}(x,x_0)$, for some $M>0$, since $K_0 e^{2u}\in C^1(\bar \Om)$. Thus, we assume that $u(x)<1$ and $u(x_0)\ge 1$. In this case, $|A_{\bar u}(x)-A_{\bar u}(x_0)|= |A_{\bar u}(x)|=|K_0(x)|(1-e^{2u})\le (\max|K_0|)(e^{2u(x_0)}-e^{2u(x)})\le M d_{g_0}(x,x_0)$, for some $M>0$, since $u\in C^1(\bar\Om)$.

Now, consider $\bar v:\bar \Om\to \real$ a solution of the Dirichlet problem  
\begin{equation*}
\De_{g_0} \bar v = A_{\bar u}, \mbox{ in } \Om, \mbox{ and } \bar v|_{\p\Om}=0.
\end{equation*} 
Since $A_{\bar u}$ is Lipschitz in $\bar \Om$ we have $v\in C^2(\bar \Om)$ (see Theorem 8.34, pg 211, of \cite{GT}). Furthermore, since $\De_{g_0}(\bar u-\bar v)\le 0$ and $(\bar u-\bar v)|_{\p\Om}=0$, one has  $\bar v \le \bar u \le 0$. This implies $1-e^{2\bar u}\le 1-e^{2\bar v}$ and $\chi_{\{\bar u<0\}}\le \chi_{\{\bar v<0\}}$, hence $A_{\bar u}\le A_{\bar v}$ in $\bar \Om$, since $K_{g_0}\ge 0$. Thus, 
\begin{equation*}
\De_{g_0}\bar v \le A_{\bar v},  \mbox{ in } \Om, \mbox{ and } \bar v|_{\p\Om}=0.
\end{equation*}
Let $v$ be defined by 
\begin{equation*}
v(x) = \left\{
\begin{array}{l}
\bar v(x), \mbox{ if } x\in \bar \Om;\\
0, \mbox{ if } x\in M\setminus\bar\Om.
\end{array}\right.
\end{equation*} 
We have $A_v$ is Lipschitz and $v$ and satisfies 
\begin{equation*}
\De_{g_0}v \le A_v, \mbox{ in } M,  
\end{equation*}
in the sense of distributions, and $v|_{\p M}=0$. Let $\om\in C^2(M)$ be a solution of
\begin{equation*}
\De_{g_0}\om = A_v,  \mbox{ in } M, \mbox{ and } \om|_{\p M}=0.
\end{equation*}
Since $\Om$ is a domain in $M$, it follows $v=0$ in a neighborhood $\m U$ of $\p M$ in $M$, hence $A_v=0$ in $\m U$, hence $w\in C^2(M)$. Furthermore, we have $\De_{g_0}(v-w)\le 0$, $(v-w)|_{\p M}=0$ and $v=0$ in $\m U$. These imply  $w\le v\le 0$ and 
\begin{equation}\label{null-deriv}
\frac{\p \om}{\p\eta}\geq \frac{\p v}{\p \eta}=0 \mbox{ on } \p M.
\end{equation}
In addition, we also have $A_v\le A_\om$. Hence, $\De_{g_0}\om \le A_\om$.
Thus, the metric $\ti g=e^{2\om}g_0$ satisfies 
\begin{eqnarray*}
K_{\ti g} &=& e^{-2\om}(K_{g_0}-\De_{g_0}\om)\ge e^{-2\om}(K_{g_0} - A_\om(1-e^{2\om})) 
\\&\ge& e^{-2\om}(K_{g_0} - \chi_{\{\om<0\}}K_{g_0}(1-e^{2\om}))
\\&=& e^{-2\om}(K_{g_0}(1 - \chi_{\{\om<0\}}) + \chi_{\{\om<0\}} K_{g_0}e^{2\om})
\\&=& K_{g_0}.
\end{eqnarray*}
The last equality follows just analyzing the cases $\om<0$ and $\om=0$. And, by (\ref{null-deriv}), one has $k_{\ti g} = \frac{\p \om}{\p\eta}+k_{g_0}\ge k_{g_0}$.

Since $\om\le 0$ and $\m L_{g_0}f\le 0$, for some $f\in C^2(M\setminus \p M)\cap C^1(M)$ positive a.e., we obtain $\m L_{g_0}f(1-e^{-2\om})\ge 0$. By Theorem \ref{min-oo conformal}, one has $\ti g=g_0$, hence $\om=0$. This implies  $v=\bar v=\bar u=0$, hence $u\ge 0$. We obtain $g\ge g_0$. Moreover, using $u\ge 0$ and $u|_{\p\Om}=0$, by (\ref{meancurvaturen=2}), one has $\frac{\p u}{\p \eta}\le 0$, hence $\ka_g\le \ka_{g_0}$. 

Now, assume further $K_g\ge K_{g_0}$ in $\Om$. Using (\ref{supharmonic}), one has $\De_{g_0}u\le 0$. Since $\frac{\p u}{\p\eta}\le 0$ and $u|_{\p\Om}=0$,  by interior maximum principle and Hopf Lemma, it follows $u=0$ in $\bar \Om$, provided $u=0$ somewhere in $\Om$ or $\frac{\p u}{\p\eta}=0$ somewhere on $\p\Om$.  

Now, consider the case $n\ge 3$ and write $g=u^{\frac{4}{n-2}}g_0$, for some positive function $u\in C^2(\bar\Om)$ with $u=1$ on $\p\Om$. 
The function $\bar u = \min\{1,u\}$ is continuous in $\bar\Om$ and satisfies $\bar u|_{\p\Om}=1$. Furthermore,  $\De_{g_0}\bar u \le \chi_{\{u<1\}}\De_{g_0}u$, in the sense of distributions. Thus, using $\chi_{\{u<1\}} R_{g}\ge \chi_{\{u<1\}} R_{g_0}\ge 0$, by (\ref{conformal n>2}), we obtain
\begin{eqnarray}\label{nge3Laplacian}
\De_{g_0}\bar u &\le&  \frac{n-2}{4(n-1)}\chi_{\{u<1\}} (R_{g_0}u -  R_{g} u^{\frac{n+2}{n-2}}) 
 \\&\le&   \frac{n-2}{4(n-1)}\chi_{\{\bar u<1\}} R_{g_0}(\bar u - \bar u^{\frac{n+2}{n-2}})\nonumber
\\ \nonumber&=& A_{\bar u}\,\bar u, \mbox{ in } M,  
\end{eqnarray}
in the sense of distributions, where $A_{\bar u}= \frac{n-2}{4(n-1)} \chi_{\{\bar u<1\}} R_{g_0}(1-{\bar u}^{\frac{4}{n-2}})$. Note that  $A_{\bar u}\ge 0$ is Lipschitz in $\bar \Om$. 

Let $\bar v\in C^2(\bar\Om)$ be a solution of the Dirichlet problem
\begin{equation*}
\De_{g_0}\bar v - A_{\bar u}\,\bar v =0 \ \mbox{ and }\  \bar v|_{\p\Om}=1
\end{equation*}
(see Theorem 8.34, pg 211, of \cite{GT}). Since $\De_{g_0}(\bar v-\bar u) - A_{\bar u}(\bar v-\bar u) \ge 0$ and $(\bar u-\bar v)|_{\p\Om}=0$, we have $\bar v\le \bar u$, since, by the strong maximum principle, $\bar v-\bar u$ cannot achieve a nonnegative maximum in the interior of $\Om$, unless $\bar u=\bar v$. We obtain $\chi_{\{\bar v<1\}}\ge \chi_{\{\bar u<1\}}$ and $1-\bar v^{\frac{4}{n-2}}\ge 1-\bar u^{\frac{4}{n-2}}$. This implies $A_{\bar v} \ge A_{\bar u}$.  Hence,
\begin{equation*}
\De_{g_0}\bar v - A_{\bar v}\,\bar v \le 0, \mbox{ in } M, \ \mbox{ and }\  \bar v|_{\p\Om}=1.
\end{equation*}

Let $v:M\to \real$ be defined by 
\begin{equation*}
v(x) = \left\{
\begin{array}{l}
\bar v(x), \mbox{ if } x\in \bar \Om;\\
1, \mbox{ if } x\in M\setminus\bar\Om.
\end{array}\right.
\end{equation*} 
Note that $v\le 1$ in $M$ and $A_v$ is Lipschitz.  Furthermore, it holds
\begin{equation}\label{dirichlet-v}
\De_{g_0} v - A_v v \le 0, \mbox{ in } M,
\end{equation}

Let $w \in C^2(M)$ be a solution of the Dirichlet problem
\begin{equation}\label{dirichlet-w}
\De_{g_0}w - A_v \, w=0, \mbox{ in } M, \mbox{ and } w|_{\p M}=1.
\end{equation}
Since $A_v\ge 0$, by the strong maximum principle, $-w$ cannot achieve a nonnegative maximum in $M\setminus \p M$, unless $w$ is constant. Hence $w>0$, since $w|_{\p M}=1$. Furthermore, by (\ref{dirichlet-v}) and (\ref{dirichlet-w}), we have $\De_{g_0}(w-v)-A_v(w-v)\ge 0$, in $M$, in the sense of distributions, and $w-v=0$ in $\p M$. Again by the strong maximum principle, we obtain $w\le v\le 1$ in $M$, hence $A_w\ge A_v$. Thus, by (\ref{dirichlet-w}),
\begin{equation}
\De_{g_0}w - A_w \, w \le 0, \mbox{ in } M, \mbox{ and } w|_{\p M}=1.
\end{equation}
Consider the metric $\ti g = w^{\frac{4}{n-2}}g_0$. By (\ref{conformal n>2}), the scalar curvatures $R_{\ti g}$ and $R_{g_0}$ satisfy
\begin{eqnarray*}
R_{\ti g} &=& w^{-\frac{n+2}{n-2}}(R_{g_0}w - \frac{4(n-1)}{n-2}\De_{g_0}w)
\\&\ge& w^{-\frac{n+2}{n-2}}(R_{g_0}w - \frac{4(n-1)}{n-2}A_w\, w)
\\&=& w^{-\frac{n+2}{n-2}}((1-\chi_{\{w<1\}})R_{g_0}w + \chi_{\{w<1\}}R_{g_0}w^{\frac{n+2}{n-2}})
\\&=& R_{g_0}.
\end{eqnarray*}
The last equality follows just by analyzing the cases $w<1$ and $w=1$. Furthermore, since $w\le 1$ and $w|_{\p M}=1$, we have $\frac{\p w}{\p \eta}\ge  0$ on $\p M$. By (\ref{meancurvature n>2}), the mean curvatures $H_{g_0}$ and $H_{\ti g}$ satisfy 
\begin{eqnarray*}
H_{\ti g} = H_{g_0}+ \frac{2(n-1)}{n-2}\frac{\p w}{\p \eta}\ge H_{g_0}.
\end{eqnarray*}

Since 
$1-w^{\frac{-4}{n-2}}\le 0$, and there exists $f\in C^2(M\setminus \p M)\cap C^1(M)$ positive a.e. such that $\m L_{g_0}f\le 0$, we have $\m L_{g_0}f (1-w^{\frac{-4}{n-2}})\ge 0$. By Theorem \ref{min-oo conformal}, it holds $\ti g=g_0$, hence $w=1$. Thus $v=\bar v=\bar u=1$, hence $u\ge 1$. This implies $g\ge g_0$. Moreover, since $u\ge 1$ and $u|_{\p M}=1$, we also have $\frac{\p u}{\p \eta}\le 0$, hence $H_g\le H_{g_0}$.

Now, we assume further $R_g \ge R_{g_0}$ in $M$. Since $u\ge 1$, by (\ref{conformal n>2}), one has $\De_{g_0}u\ge 0$. Thus, if $u=1$, somewhere in $\Om$, or $\frac{\p u}{\p \eta}=0$, somewhere in $\p\Om$, then, by interior maximum principle or Hopf Lemma, it holds $u=1$ in $M$. Theorem \ref{thm-hang-wang} is proved.

\end{document}